\documentclass[a4paper, 11pt]{amsart}
\usepackage[latin1]{inputenc}
\usepackage[ T1]{fontenc}
\usepackage[english]{babel}
\usepackage{amssymb}
\usepackage{amsmath}
\usepackage{amsthm}
\usepackage{amscd}
\usepackage{amsfonts}
\usepackage{stmaryrd}
\usepackage{pb-diagram}
\usepackage{epic,eepic,epsfig}
\usepackage{a4wide}
\usepackage{nextpage}
\usepackage{fancyhdr}
\usepackage{enumerate}
\usepackage{hyperref}
\usepackage{float}

\newtheorem{prop}{Proposition}[section]
\newtheorem{corollary}[prop]{Corollary}

\newtheorem{lemma}[prop]{Lemma}
\newtheorem{thm}[prop]{Theorem}

\theoremstyle{definition}
\newtheorem{rem}[prop]{Remark}

\newcommand{\SL} {\mathrm{SL}}
\newcommand{\GL} {\mathop{\mathrm{GL}}}

\renewcommand{\Im}{\mathop{\mathrm{Im}}}

\def\BC{\mathbb{C}}
\def\BF{\mathbb{F}}
\def\BQ{\mathbb{Q}}
\def\BZ{\mathbb{Z}}
\def\BH{\mathbb{H}}

\def\BT{\mathcal{BT}}
\def\CF{\mathcal{F}}

\def\Cinf{\mathbb{C}_\infty}
\def\diag{\mathrm{diag}}
\def\GL{\mathrm{GL}}
\def\lcm{\mathrm{lcm}}

\usepackage[OT2,T1]{fontenc}
\DeclareSymbolFont{cyrletters}{OT2}{wncyr}{m}{n}
\DeclareMathSymbol{\Sha}{\mathalpha}{cyrletters}{"58}

\begin{document}

\title{Explicit bounds on the coefficients of modular polynomials for the elliptic $j$-invariant}
\author{}
\author{Florian Breuer \and Fabien Pazuki}

\address{School of Information and Physical Sciences, The University of Newcastle,
University Drive, Callaghan, NSW 2308, Australia.}
\email{Florian.Breuer@newcastle.edu.au}

\address{Department of Mathematical Sciences, University of Copenhagen,
Universitetsparken 5, 
2100 Copenhagen \O, Denmark, and Universit\'e de Bordeaux, 33405 Talence, France.}
\email{fpazuki@math.ku.dk}

\thanks{The authors thank the IRN GandA (CNRS). The second author is supported by ANR-20-CE40-0003 Jinvariant.}
\maketitle

\noindent \textbf{Abstract.}
We obtain an explicit upper bound on the size of the coefficients of the elliptic modular polynomials $\Phi_N$ for any $N\geq1$.
These polynomials vanish at pairs of $j$-invariants of elliptic curves linked by cyclic isogenies of degree $N$.
The main term in the bound is asymptotically optimal as $N$ tends to infinity.

{\flushleft
\textbf{Keywords:} Modular polynomials, elliptic curves.\\
\textbf{Mathematics Subject Classification:} 11G05. }

\begin{center}
---------
\end{center}

\thispagestyle{empty}

\maketitle

\section{Introduction}\label{section def}


For any non-zero polynomial $P$ in one or more variables and complex coefficients we define its {\em height} to be 
\[
h(P) := \log \max|c|, \quad\text{where $c$ ranges over all coefficients of $P$.}
\]

Let $N$ be a positive integer and denote by $\Phi_N = \Phi_N(X,Y) \in \BZ[X,Y]$
the (classical) modular polynomial, which vanishes at pairs of $j$-invariants of elliptic curves linked by a cyclic $N$-isogeny, see \cite[Chapter 5]{Lang}. 
Alternatively, if we view $j$ as the function on the complex upper half-plane where $j(\tau)$
is the $j$-invariant of the complex elliptic curve $\BC/(\BZ+\tau\BZ)$, then
$\Phi_N(X, j(\tau))$ is the minimal polynomial of $j(N\tau)$ over $\BC(j(\tau))$.

Modular polynomials have important applications in cryptography and certain algorithms for computing $\Phi_N$ require explicit bounds on the size of the coefficients, so one is interested in explicit bounds on $h(\Phi_N)$. 

Paula Cohen Tretkoff \cite{Coh} proved that when $N$ tends to $+\infty$
\begin{equation}\label{Paula}
   h(\Phi_N) = 6\psi(N)\big[\log N - 2\kappa_N + O(1)\big] 
\end{equation}
where
\[
\psi(N) = N\prod_{p|N}\left(1 + \frac{1}{p}\right)\quad\text{and}\quad
\kappa_N = \sum_{p|N}\frac{\log p}{p},
\]
but the implied bounded function is not explicit.

In the case where $N=l$ is prime, Br\"oker and Sutherland \cite{BrSu} estimated the constants in Cohen's argument to obtain
\[
h(\Phi_l) \leq 6l\log l + 16 l + 14\sqrt{l}\log l.
\]

In the general case, the second author \cite{Paz} obtained in his Corollary 4.3, via a different method,
\begin{equation}\label{Pazweak}
  h(\Phi_N) \leq \psi(N)\big[ 6\log N + \log\psi(N) + 6\log(12 \log N + 2\log\psi(N) + 25.2) + 15.7\big].  
\end{equation}

Inequality (\ref{Pazweak}) has the merit of being completely explicit for all $N\geq1$, but the 
main term is slightly too big when compared with the asymptotic of (\ref{Paula}).

The goal of the present paper is to prove the following result, where we solve this issue and provide an upper bound with the correct main term for all $N$. Let us first define 
\[
\lambda_N := \sum_{p^n\|N}\frac{p^n-1}{p^{n-1}(p^2-1)}\log p.
\]

\begin{thm}\label{Main1}
Let $N\geq 2$. The height of the modular polynomial $\Phi_N(X,Y)$ is bounded by
\begin{equation}
    \label{eq:main1}
    h(\Phi_N) \leq 6\psi(N)\big[\log N -2\lambda_N + \log\log N + 4.436\big].
\end{equation}
\end{thm}




We prove this theorem using a different path than the one followed in \cite{Paz}. The main new ingredient is a finer estimate of the Mahler measure of $j$-invariants, coming from previous work of Pascal Autissier \cite{Aut}. We also use precise analytic estimates for the discriminant modular form on the fundamental domain of the upper half plane (under the classical action of $\mathrm{SL}_2(\mathbb{Z})$), and a classical interpolation method to help us derive bounds on the height of a polynomial in two variables, from knowledge of the height of several specializations of this polynomial.

Let us now discuss the optimality of the bound. The main term is the expected one. For lower order terms, notice that
\[
-0.385 < -\sum_{p|N}\frac{\log p}{p(p+1)} \leq \lambda_N - \kappa_N \leq 
\sum_{p|N}\frac{\log p}{p(p^2-1)} < 0.186,
\]
so one changes little replacing $\lambda_N$ by $\kappa_N$ in Theorem \ref{Main1}. 
On the other hand, one would like to get rid of the spurious $\log\log N$ term, but
for practical purposes this might be less useful than keeping the constant as small as possible.

\begin{figure}
    \centering
    \includegraphics[width=\textwidth]{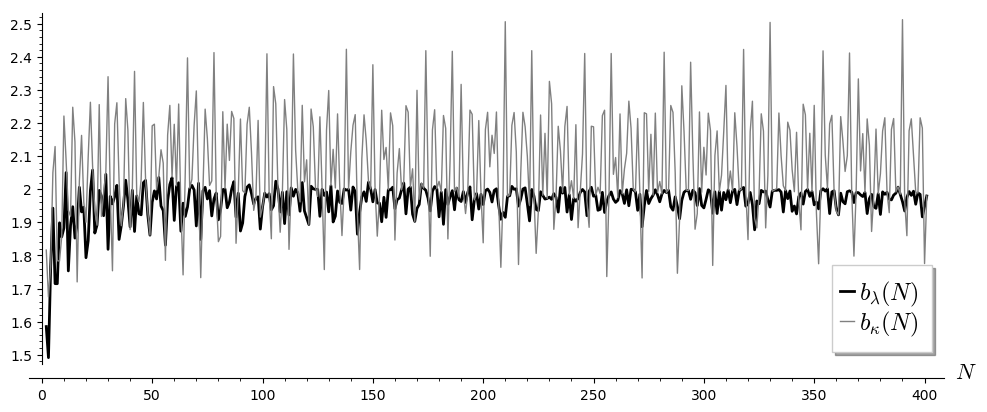}
    \caption{The bounded functions $b_\lambda(N)$ (bold) and $b_\kappa(N)$ (grey) satisfying
    $h(\Phi_N) = 6\psi(N)\big[\log N -2\lambda_N + b_\lambda(N)\big] 
    = 6\psi(N)\big[\log N -2\kappa_N + b_\kappa(N)\big]$ 
    for $N\leq 400.$ 
    Notice that $b_\lambda(N) < 2.1$ in this range. Theorem \ref{Main1} is equivalent to 
    $b_\lambda(N) \leq \log\log N + 4.436$.
    }
    \label{fig:b}
\end{figure}

It is interesting to consider the functions $b_\lambda(N)$ and $b_\kappa(N)$ for which
\begin{equation}
h(\Phi_N) = 6\psi(N)\big[\log N -2\lambda_N + b_\lambda(N)\big] =
6\psi(N)\big[\log N -2\kappa_N + b_\kappa(N)\big].
\end{equation}
These functions are plotted in Figure \ref{fig:b} for $N\leq 400$, based on computations of $\Phi_N$ by
Andrew Sutherland \cite{SuWeb} using 
the algorithms in \cite{BKS} (for prime $N$) and \cite{BOS} (for composite~$N$).

The content of Cohen's Theorem is that 
$b_\kappa(N)$ and thus also $b_\lambda(N)$ are bounded functions.
Our Theorem \ref{Main1} is equivalent to
$b_\lambda(N) \leq \log\log N + 4.436$, which is clearly seen to hold for $N\leq 400$; in fact, 
$b_\lambda(N) < 2.1$ in this range.

In our proof of Theorem \ref{Main1} we may thus assume that $N > 400$.
We explain in Remark \ref{N_0} and in Lemma \ref{lem:limit} that more computations for $N > 400$ lead to minor improvements on the constant $4.436.$

From Figure \ref{fig:b} it appears that $b_\lambda(N)$ is bounded more tightly than $b_\kappa(N)$, thus suggesting that $\lambda_N$ is a more natural function to use in the bound for $h(\Phi_N)$ than is $\kappa_N$.



\medskip
\paragraph{\bf Acknowledgements.} The authors are grateful to Pascal Autissier for suggesting that the results in \cite[\S2]{Aut} might be fruitfully applied to estimating $h(\Phi_N)$. They are also grateful to Joseph Silverman for an interesting discussion around \cite{Sil}. The authors warmly thank Andrew Sutherland for the computations of modular polynomials he performed in record time to help them improve numerical values in the statement of Theorem \ref{Main1}. They also thank the referees for very efficient feedback.   
The authors thank the IRN GandA (CNRS). The second author is supported by ANR-20-CE40-0003 Jinvariant.  

\section{Preliminary results}

Denote the complex upper half-plane by
\[
\BH := \{z\in\BC \;|\; \Im(z) > 0 \}.
\]
Every $\tau\in\BH$ defines a lattice $\Lambda_\tau = \BZ + \tau\BZ$ in $\BC$, and it is well known that
every complex elliptic curve is isomorphic to $\BC/\Lambda_\tau$ for some $\tau\in\BH$. 
If we denote the $j$-invariant of this elliptic curve by $j(\tau)$, then 
\[
j : \BH \longrightarrow \BC
\]
defines an analytic function on $\BH$.

The group $\SL_2(\BZ)$ acts on the upper half-plane $\BH$ by 
\[
\gamma(\tau) := \frac{a\tau + b}{c\tau + d}, \quad\text{where}\quad \gamma = 
\left(\begin{matrix}a & b \\ c & d \end{matrix}\right) \in \SL_2(\BZ).
\]
A fundamental domain for this action is given by
\[
\mathcal{F}=\{\tau\in\BH\; : \; \vert \tau\vert\geq 1, \;  -\frac{1}{2} < \mathrm{Re}(\tau) \leq \frac{1}{2} \; \mathrm{and}\;
\mathrm{Re}(\tau) \geq 0 \; \mathrm{if}\; |\tau|=1\}.
\]
Thus every $\tau\in\BH$ is $\SL_2(\BZ)$-equivalent to an element $\tilde{\tau}\in\CF$, which we call \emph{reduced}.

The modular function $j : \BH \rightarrow \BC$ is $\SL_2(\BZ)$-invariant. We define 
\[
q = e^{2\pi i \tau}, \quad \tau\in\BH,
\]
then the Fourier expansion at infinity of $j$ can be written as a $q$-expansion 
\[
j(\tau) = \frac{1}{q}+744+196884q+\ldots.
\]




We denote by $\Delta$ the modular discriminant function
\[
\Delta : \BH \longrightarrow \BC,
\]
which is a weight 12 cusp form for $\SL_2(\BZ)$. 
We normalize $\Delta$ so that its $q$-expansion is
\[
\Delta(\tau) = q\prod_{n=1}^{\infty}(1-q^n)^{24} = q -24q^2 + 252q^3 + \cdots.
\]
We point out that the discriminant of the elliptic curve $E_\tau$ is given by $(2\pi)^{12}\Delta(\tau)$, which is why
most sources (e.g. \cite{Lang}) normalize $\Delta$ differently, multiplying the above product by the factor $(2\pi)^{12}$.
We choose our normalization to be consistent with \cite{Paz}, which contains estimates that we will use.


Let us denote, for $N\geq1$, 
\[
C_N=\left\{\left(\begin{matrix}
a & b \\
0 & d
\end{matrix}\right)\; : \; a,b,d\in\BZ, \;  ad=N, \; a\geq1, \; 0\leq b\leq d-1,\; \gcd(a,b,d)=1 \right\}.
\]
We have
\[
\#C_N = \psi(N) = N\prod_{p|N}\left(1 + \frac{1}{p}\right).
\]
The elements of $C_N$ encode cyclic $N$-isogenies in the following way.
Let $E_\tau$ be an elliptic curve. For each
\[
\gamma = \left(\begin{matrix}
a_\gamma & b_\gamma \\
0 & d_\gamma
\end{matrix}\right)\in C_N,
\]
we let 
\[
\tau_\gamma = \gamma(\tau) = \frac{a_\gamma\tau + b_\gamma}{d_\gamma}, \quad
\Lambda_\gamma = \BZ + \tau_\gamma\BZ, \quad\text{and}\quad 
E_\gamma = E_{\tau_\gamma} = \BC/\Lambda_\gamma.
\]
Then the natural map
\[
E \longrightarrow E_\gamma, \quad (z \bmod \Lambda_\tau) \longmapsto (z \bmod \Lambda_\gamma)
\]
is a cyclic $N$-isogeny. 

Furthermore, up to isomorphism, every cyclic $N$-isogeny with source $E$ arises in this way.
In particular, we have the factorization
\[
\Phi_N\big(X, j(\tau)\big)=\prod_{\gamma\in{C_N}}\big(X-j(\tau_\gamma)\big).
\]


Our goal is to 
bound the coefficients 
of the modular polynomial $\Phi_N(X,Y)$. By interpolation, it is enough to estimate the height of $\Phi_N(X, j(\tau))$ for several carefully chosen $\tau\in\BH$. 

By \cite[Lemma 1.6]{BrZu}
the height of $\Phi_N(X,j(\tau))$ is bounded in terms of its Mahler measure 
\begin{equation}\label{SN}
    S_N(\tau)=\displaystyle{\sum_{\gamma\in{C_N}}\log\max\big(1, \vert j(\tau_\gamma)\vert\big)}
\end{equation}
by 
\begin{equation}
    \label{eq:mahler}
    h(\Phi_N(X,j)) \leq S_N(\tau) + \log \binom{\psi(N)}{\psi(N)/2}
    \leq S_N(\tau) + \psi(N)\log 2.
\end{equation}

We will concentrate on estimating $S_N(\tau)$ for a fixed $\tau\in\BH$.


In general, $\tau_\gamma$ won't be reduced, so we choose 
$\left(\begin{matrix}a & b \\ c & d\end{matrix}\right)\in\SL_2(\BZ)$
for which 
\[
\tilde{\tau}_\gamma = \frac{a\tau_\gamma + b}{c\tau_\gamma +d} \in \CF
\]
is reduced. 
%
%
Since 
\[
\Im(\tilde{\tau}_\gamma) = \Im\left(\frac{a\tau_\gamma + b}{c\tau_\gamma + d}\right) = 
\frac{\Im(\tau_\gamma)}{|c\tau_\gamma + d|^2},
\]
we obtain 
\begin{equation} \label{eq:cIm}
    -\log|c\tau_\gamma + d| = \frac{1}{2}\big[\log\Im(\tilde{\tau}_\gamma) - \log\Im(\tau_\gamma)\big].
\end{equation}

Also, since $\Delta$ is a modular form of weight $12$ for $\SL_2(\BZ)$, we find that 
\[
\tilde{\Delta}_\gamma := \Delta(\tilde\tau_\gamma) = (c\tau_\gamma + d)^{12}\Delta(\tau_\gamma) =: 
(c\tau_\gamma + d)^{12}\Delta_\gamma,
\]

so
\begin{align} 
   \log|\Delta_\gamma| &= \log|\tilde{\Delta}_\gamma| - 12\log|c\tau_\gamma + d| \nonumber \\
   &= \log|\tilde{\Delta}_\gamma|+ 6\big[ \log\Im(\tilde{\tau}_\gamma) - \log\Im(\tau_\gamma) \big]. \label{eq:cDelta}
\end{align}

Note that by \cite[Lemma 2.4]{Paz} we have
\begin{equation}\label{eq:AIL}
    \log\Im(\tilde{\tau}_\gamma) - \log\Im(\tau) \leq \log N
\end{equation}
for each $\gamma\in C_N$, provided that $\tau\in\CF$.

We need a few more preliminaries:

By \cite[Lemme 2.2]{Aut}, we have
\[
\prod_{\gamma\in C_N} \Delta(\gamma(\tau)) = \big[-\Delta(\tau)\big]^{\psi(N)},
\]
so we get
\begin{equation}
    \label{eq:Aut2.2}
    \sum_{\gamma\in C_N} \log|\Delta_\gamma| = \psi(N)\log|\Delta|.
\end{equation}

Furthermore, \cite[Lemme 2.3]{Aut} says
\[
  \sum_{\gamma\in C_N}\log\frac{d_\gamma}{a_\gamma} = \psi(N)(\log N - 2\lambda_N),
\]
which combined with 
\[
\Im(\tau_\gamma) = \Im\left( \frac{a_\gamma\tau + b_\gamma}{d_\gamma} \right)
= \frac{a_\gamma}{d_\gamma}\Im(\tau)
\]
gives
\begin{equation} 
  \label{eq:Aut2.3}
  -\sum_{\gamma\in C_N}\log\Im(\tau_\gamma) 
  = \psi(N)\big(\log N - 2\lambda_N   -\log\Im(\tau)\big).
\end{equation}

Finally, since $\tilde{\tau}_\gamma\in\CF$, \cite[(2.22)]{Paz} gives us, if we denote $j_\gamma=j(\tau_\gamma)$,
\begin{equation}
   \label{eq:Paz2.22}
    \Im(\tilde{\tau}_\gamma) \leq \frac{1}{2\pi}\log(|j_\gamma|+970.8),
\end{equation}
whereas \cite[(3.18)]{Paz} gives, for any $\gamma\in C_N$,
\begin{equation}
    \label{eq:Paz3.18}
    \log\max(|\tilde{\Delta}_\gamma|,|j_\gamma\tilde{\Delta}_\gamma|) \leq \log(9.02).
\end{equation}
This last estimate depends on our choice of normalisation of $\Delta(\tau)$.

We note that the identities (\ref{eq:Aut2.2}) and (\ref{eq:Aut2.3}) from \cite{Aut} involve the non-reduced $\tau_\gamma$, 
whereas the estimates (\ref{eq:AIL}), (\ref{eq:Paz2.22}) and (\ref{eq:Paz3.18}) from \cite{Paz} depend on the reduced $\tilde{\tau}_\gamma$. 
The main idea of this paper is to combine these ingredients using~(\ref{eq:cDelta}).

\section{Proof of Theorem \ref{Main1}}

We are now ready to start our main calculation on the sum $S_N(\tau)$ from (\ref{SN}).
\begin{align}
  S_N(\tau)= & \sum_{\gamma\in C_N}\log\max(|\Delta_\gamma|,|j_\gamma\Delta_\gamma|)
      - \sum_{\gamma\in C_N}\log|\Delta_\gamma| \nonumber\\
    = & \sum_{\gamma\in C_N}\log\max(|\Delta_\gamma|,|j_\gamma\Delta_\gamma|)
      - \psi(N)\log|\Delta|  \quad \text{(by (\ref{eq:Aut2.2}))}  \nonumber\\
    = & \sum_{\gamma\in C_N} \log\max(|\tilde{\Delta}_\gamma|,|j_\gamma\tilde{\Delta}_\gamma|) 
    +6 \sum_{\gamma\in C_N} \big[\log\Im(\tilde{\tau}_\gamma) 
    - \log\Im(\tau_\gamma)\big]
    - \psi(N)\log|\Delta|  \quad \text{(by (\ref{eq:cDelta})),} \nonumber\\
    \end{align}  
    hence we get
    \begin{align}
    S_N(\tau)\leq\; & \psi(N)\log(9.02)
    +6 \sum_{\gamma\in C_N} \big[\log\Im(\tilde{\tau}_\gamma) 
    - \log\Im(\tau_\gamma)\big]
    - \psi(N)\log|\Delta|  \quad \text{(by (\ref{eq:Paz3.18}))} \nonumber\\
    = \; & \psi(N)\log(9.02) 
    + 6\psi(N)\big(\log N - 2\lambda_N   -\log\Im\tau\big) \nonumber\\
    & + 6\sum_{\gamma\in C_N} \log\Im(\tilde{\tau}_\gamma) 
    - \psi(N)\log|\Delta| \quad\text{(by (\ref{eq:Aut2.3}))}
    \nonumber\\
   \leq \; & 6\psi(N)\big[\log N - 2\lambda_N +  0.367\big]
    + 6\sum_{\gamma\in C_N} \log\Im(\tilde{\tau}_\gamma) 
    - \psi(N)\log\big[|\Delta|(\Im\tau)^6\big].
    \label{eq:halfway}
\end{align}    

At this point we record the following intermediate result.
If $\tau\in\CF$ then we may apply (\ref{eq:AIL}) and obtain 
\begin{align}
    S_N(\tau) 
    \leq\; & 
    6\psi(N)\big[\log N - 2\lambda_N +  0.367\big] 
    +6\psi(N)[\log N + \log\Im\tau] - \psi(N)\log\big[|\Delta|(\Im\tau)^6\big] \nonumber\\
     \leq\; & \psi(N)[12\log N + 2.199 - \log|\Delta|]. 
    \label{eq:crude}
\end{align}


We continue our calculation from (\ref{eq:halfway}).
    
\begin{align*}
    S_N(\tau) 
    \leq\; & 
    6\psi(N)\big[ \log N -2\lambda_N + 0.367\big]
    - \psi(N)\log\big[|\Delta|(\Im\tau)^6\big] \\
    & +6 \sum_{\gamma\in C_N} \log\Big[\frac{1}{2\pi}\log(|j_\gamma| + 970.8)\Big] 
     \quad \text{(by (\ref{eq:Paz2.22}))} \\
    =\; & 6\psi(N)\big[\log N - 2\lambda_N + 0.367\big] 
    -\psi(N)\log\big[|\Delta|\Im(\tau)^6\big] \\
    & + 6\psi(N)\log\prod_{\gamma\in C_N} \Big[\frac{1}{2\pi}\log(|j_\gamma| + 970.8)\Big]^{1/\psi(N)} \\
    \leq\; &  6\psi(N)\big[\log N - 2\lambda_N + 0.367\big] 
    -\psi(N)\log\big[|\Delta|\Im(\tau)^6\big] \\
    & + 6\psi(N)\log \Big[ \frac{1}{2\pi\psi(N)} \sum_{\gamma\in C_N} \log(|j_\gamma| + 970.8) \Big],  
\end{align*}
where the last inequality follows by the arithmetic-geometric mean inequality.

For any real number $x$, the inequality $x+970.8\leq971.8\max\{1,x\}$ holds, so we finally obtain 
\begin{align}
    S_N(\tau) =\; & \sum_{\gamma\in C_N}\log\max(1,|j_\gamma|) \nonumber \\
    \leq\; &   6\psi(N)\big[\log N - 2\lambda_N + 0.367\big] 
    -\psi(N)\log\big[|\Delta|\Im(\tau)^6\big] \nonumber \\
    & + 6\psi(N)\big[ \log S_N(\tau) + \log\log(971.8) - \log\psi(N) - \log(2\pi)\big] \nonumber \\
    \leq\; & 6\psi(N)\big[\log N - 2\lambda_N + \log\big(S_N(\tau)/\psi(N)\big) + 0.458]
    -\psi(N)\log\big[|\Delta|\Im(\tau)^6\big]. \label{eq:semifinal}
\end{align}


To deduce an explicit bound on $S_N(\tau)$, we start with a crude bound on $S_N(\tau)/\psi(N)$, then strengthen our result recursively. More precisely, we prove the following technical lemma.

\begin{lemma}\label{lem:induction}
    Fix $\tau\in\BH$ and let
    \begin{align*}
        a(\tau) &= 0.458 - \frac{1}{6}\log\big[|\Delta(\tau)|\Im(\tau)^6\big] \\
        b(\tau) &= 2.199 - \log|\Delta(\tau)|.
    \end{align*}

    Suppose that $N > N_0 \geq 3$. 
    Consider the sequence $\big(c_n(\tau)\big)_{n\geq 0}$ defined recursively by
    \begin{align*}
        c_0(\tau) =&\; a(\tau) + \log\left[12 + \frac{b(\tau)}{\log N_0}\right], \\
        c_{n+1}(\tau) =&\; a(\tau) + \log 6 + 
        \log\left[1 + \frac{\log\log N_0 + c_n(\tau)}{\log N_0}\right], \quad n\geq 0.
    \end{align*}
    Then for all $n\geq 0$,
    \begin{equation}
        \label{eq:B}
        S_N(\tau) \leq 6\psi(N)\big[\log N - 2\lambda_N + \log\log N + c_n(\tau) \big].
    \end{equation}
\end{lemma}

\begin{proof}
The bound (\ref{eq:crude}) gives 
\begin{align*}
S_N(\tau)/\psi(N) 
&\leq 12\log N + b(\tau) \\
&\leq \left[12 + \frac{b(\tau)}{\log N_0}\right]\log N.
\end{align*}
Plugging this into (\ref{eq:semifinal}) gives us (\ref{eq:B}) with $n=0$.

Next, assume (\ref{eq:B}) holds for some $n\geq 0$. Since $N > N_0$, we obtain
\[
\log\log N + c_n(\tau) < \left(\frac{\log\log N_0 + c_n(\tau)}{\log N_0}\right) \log N,
\]
so (\ref{eq:B}) gives us
\[
    S_N(\tau) \leq 6\psi(N)\left[\log N + 
    \left(\frac{\log\log N_0 + c_n(\tau)}{\log N_0}\right)\log N\right] 
\]
and so
\[
S_N(\tau)/\psi(N) \leq 6\left[1 + \frac{\log\log N_0 + c_n(\tau)}{\log N_0}\right]\log N.
\]
Plugging this into (\ref{eq:semifinal}) gives us
\[
S_N(\tau) \leq 6\psi(N)\big[\log N - 2\lambda_N + \log\log N + c_{n+1}(\tau) \big].
\]
\end{proof}

The interpolation lemma \cite[Lemma 20]{BrSu} gives, for real $L > 1$,
\[
h(\Phi_N(X,Y)) \leq \max_{L \leq j \leq 2L}h(\Phi_N(X,j)) + \psi(N)\left(\frac{\log L + 1}{L} + 3\log 2\right),
\]
so by (\ref{eq:mahler}) we get
\begin{equation}\label{interpolation}
h(\Phi_N(X,Y)) \leq \max_{L \leq j(\tau) \leq 2L} S_N(\tau) + \psi(N)\left(\frac{\log L + 1}{L} + 4\log 2\right).
\end{equation}

It is well-known that the $j$-function takes non-negative real values on the following path on the boundary of the fundamental domain $\CF$:
\[
\Gamma := \{e^{i\theta} \; | \; \frac{\pi}{3} \leq \theta \leq \frac{\pi}{2}\} \cup
\{ix \;|\; x\in [0,\infty)\}
\]
and the function $j:\Gamma \rightarrow [0,\infty)$ is a bijection.

We now define, for the values $c_n(\tau)$ in Lemma \ref{lem:induction} with $N_0=400$, 
\[
c(\tau) := \inf_{n\geq 0}c_n(\tau).
\]


Optimizing on the interval $L\leq j\leq 2L$, we obtain

\[
h(\Phi_N(X,Y)) \leq 
6\psi(N)\big[\log N - 2\lambda_N + \log\log N + c_n(\tau) \big]\big|_{j(\tau)=2L}
+ \psi(N)\left(\frac{\log L + 1}{L} + 4\log 2\right).
\]

Optimizing $c(\tau)$ (using SageMath \cite{Sage}) when $L > 1$, we obtain the strongest upper bound when we choose $L=166.48$,
then $\tau = j^{-1}(L) = e^{i\cdot 1.257}$ and
\[
a(\tau) \leq 1.5004, \quad b(\tau) \leq 8.1532, \quad c(\tau) \leq  3.9655.
\]
Putting all of this together, we obtain 
\[
h(\Phi_N) \leq 6\psi(N)\big[\log N -2\lambda_N + \log\log N + 4.436\big]
\]
for $N\geq N_0=400$. 


As can be seen from Figure \ref{fig:b}, the result also holds for $N\leq 400$, thus completing the proof of Theorem \ref{Main1}.

\qed

Let us add the following remark.

\begin{rem}\label{N_0}
The constant in Theorem \ref{Main1} can be further improved if we assume $N > N_0$ for larger values of $N_0$ and check the result for $N \leq N_0$ via direct computation.
We list below the values of the constant in Theorem \ref{Main1} obtained assuming $N > N_0$ for some other values of $N_0$.

\begin{center}
\begin{tabular}{cc}
    $N_0$: & constant: \\
    $400$ & $4.436$ \\
    $500$ & $4.418$ \\
    $1000$ & $4.373$ \\
    $2000$ & $4.336$ \\
    $5000$ & $4.292$
\end{tabular}
\end{center}

The best value is always obtained when $L=166.48$. The gain is somehow limited, even asymptotically, as explained in the next lemma. The next inequality is weaker numerically, but helps understand how the estimates on $c_n$ will evolve when $n\to+\infty$ and $N_0\to+\infty$.

\end{rem}


\begin{lemma}\label{lem:limit}
    Suppose that $N_0 \geq 3$. 
    Consider the sequence $\big(c_n(\tau)\big)_{n\geq 0}$ defined in Lemma \ref{lem:induction}.
    Then for all $n\geq 0$,
    \begin{equation}
        \label{eq:C}
        c_n(\tau)\leq \frac{c_0(\tau)}{(\log N_0)^n} + (a(\tau)+\log 6)\frac{\log N_0}{\log N_0 - 1}+\frac{\log\log N_0}{\log N_0-1}.
    \end{equation}
\end{lemma}

\begin{proof}
Let us denote $A=a(\tau)+\log 6$, for any $x\geq0$, we have $\log(1+x)\leq x$, hence we get $$c_{n+1}(\tau)\leq A+\frac{\log\log N_0+c_n(\tau)}{\log N_0},$$ which gives by induction $$c_n(\tau)\leq \frac{c_0(\tau)}{(\log N_0)^n}+\left(A+\frac{\log\log N_0}{\log N_0}\right)\sum_{k=0}^{n}\frac{1}{(\log N_0)^k}\leq \frac{c_0(\tau)}{(\log N_0)^n}+\left(A+\frac{\log\log N_0}{\log N_0}\right)\frac{\log N_0}{\log N_0 - 1},$$ which gives the conclusion.
\end{proof}

If one takes $\tau = e^{i\cdot 1.257}$ and $n\to+\infty$ in (\ref{eq:C}) we obtain the following inequality, valid for any $N_0\geq3$ and any $N\geq N_0$:

\begin{equation}\label{last}
    h(\Phi_N) \leq 6\psi(N)\left[\log N - 2\lambda_N + \log\log N + 3.293\frac{\log N_0}{\log N_0 -1} + \frac{\log\log N_0}{\log N_0 -1} +0.46537\right]. 
\end{equation}

Explicit computation of $c_n(\tau)$ will generally give better numerical values of course, but this equation (\ref{last}) gives an idea of how these estimates will vary with $N_0$.

\end{document}